\documentstyle[11pt,figure]{article}
\figon

\begin{document}

\font\bbbld=msbm10 scaled\magstep1
\newcommand{\bfR}{\hbox{\bbbld R}}
\newcommand{\bfC}{\hbox{\bbbld C}}
\newcommand{\bfZ}{\hbox{\bbbld Z}}
\newcommand{\bfH}{\hbox{\bbbld H}}
\newcommand{\bfQ}{\hbox{\bbbld Q}}
\newcommand{\bfN}{\hbox{\bbbld N}}
\newcommand{\bfP}{\hbox{\bbbld P}}
\newcommand{\bfT}{\hbox{\bbbld T}}
\def\Sym{\mathop{\rm Sym}}
\newcommand{\suchthat}{\mid}
\newcommand{\halo}[1]{\Int(#1)}
\def\Int{\mathop{\rm Int}}
\def\Re{\mathop{\rm Re}}
\def\Im{\mathop{\rm Im}}
\newcommand{\union}{\cup}
\newcommand{\goesto}{\rightarrow}
\newcommand{\bdy}{\partial}
\newcommand{\n}{\noindent}
\newcommand{\p}{\hspace*{\parindent}}

\newtheorem{theorem}{Theorem}[section]
\newtheorem{assertion}{Assertion}[section]
\newtheorem{proposition}{Proposition}[section]

\newtheorem{lemma}{Lemma}[section]
\newtheorem{definition}{Definition}[section]
\newtheorem{claim}{Claim}[section]
\newtheorem{corollary}{Corollary}[section]
\newtheorem{observation}{Observation}[section]
\newtheorem{conjecture}{Conjecture}[section]
\newtheorem{question}{Question}[section]
\newtheorem{example}{Example}[section]

\newbox\qedbox
\setbox\qedbox=\hbox{$\Box$}
\newenvironment{proof}{\smallskip\noindent{\bf Proof.}\hskip \labelsep}%
			{\hfill\penalty10000\copy\qedbox\par\medskip}
\newenvironment{remark}{\smallskip\noindent{\bf Remark.}\hskip \labelsep}%
			{\hfill\penalty10000\copy\qedbox\par\medskip}
\newenvironment{proofspec}[1]%
		      {\smallskip\noindent{\bf Proof of #1.}\hskip \labelsep}%
			{\nobreak\hfill\hfill\nobreak\copy\qedbox\par\medskip}
\newenvironment{acknowledgements}{\smallskip\noindent{\bf Acknowledgements.}%
	\hskip\labelsep}{}
\setlength{\baselineskip}{1.2\baselineskip}
\title{Limit Surfaces of Riemann Examples}
\author{David Hoffman, Wayne Rossman}
\date{}

\maketitle


\section{Introduction}

The only connected minimal surfaces foliated by circles and lines are 
domains on one 
of the following surfaces: the helicoid, the catenoid, the plane, and the 
examples of Riemann (\cite{Ri} p329-33, \cite{En} p403-6, 
\cite{Ni} p85-6).  All these surfaces are complete and embedded.  
Topologically they are planar domains: the helicoid is simply-connected, 
the catenoid is an annulus (conformally a twice-punctured sphere), and each 
Riemann example (see Figure 1) is conformal to the plane minus the points 
$\{ (n,0), \, (\frac{1}{n},0) \, | \, n \in \bfZ \}$ \cite{HKR}.  
In this section we will 
show that the plane, helicoid, and catenoid arise naturally as limits of 
well-chosen and properly normalized sequences of Riemann examples.  
The local behavior of domains on Riemann examples accounts for 
the existence of these limits, thus allowing a change of topology in 
the limit surfaces.  

\begin{theorem}
	Any sequence of Riemann examples that converges to a surface in 
$\bfR^3$ converges either to a Riemann example, a helicoid, a catenoid, 
an infinite set of equally spaced parallel planes, or a single plane.  
Furthermore, there exist sequences of Riemann examples that converge 
to each of these possibilities.  
\end{theorem}
\begin{figure}[t]
\hspace{1.3in}
\epsfysize = 2.5in	\epsffile{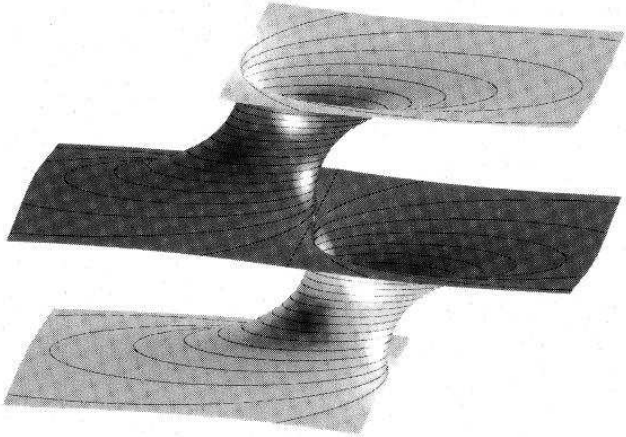}
\hfill \break
\caption{The Riemann Example ${\cal R}_\lambda, \; \lambda = 1$}
\end{figure}
In order to be precise, we state what we mean by convergence of 
surfaces in $\bfR^3$ in this context.  

\begin{definition}
	A sequence of surfaces $\{{\cal S}_i\}_{i=1}^\infty$ 
{\em converges} as $i \rightarrow \infty$ to a surface $\cal S$ in 
$\bfR^3$ if, 
for any compact region $B \subseteq \bfR^3$, there exists an 
integer $N_B$ such that for $i > N_B$, ${\cal S}_i \cap B$ 
is a normal graph over 
$\cal S$, and $\{ {\cal S}_i \cap B \}_{i=N_B}^\infty$ 
converges to $\cal S$ in the $C^\infty$-topology.  
\end{definition}

One can define analogously the convergence of a one-parameter family 
${\cal S}_t, t \in \bfR$, as 
$t \rightarrow t_0$.  Since this convergence is in the 
$C^\infty$-topology, it requires that derivatives of all orders 
of the coordinate functions on 
the graphs converge to the corresponding derivatives of the coordinate 
functions on the limit-surface.  

We observe that catenoids, helicoids, and planes are 
limit-surfaces of sequences of Riemann examples.  Therefore the 
collection of connected embedded minimal surfaces foliated 
by circles and lines is itself connected in the topology associated with 
this convergence.  

We wish to thank Rob Kusner, Pascal Romon, Ed Thayer, Johannes Nitsche
and Harold Rosenberg for helpful 
conversations.  

\section{The Family of Riemann Examples} 

\subsection{Definition of the Riemann Examples} 

	We now describe properties of the Riemann examples, and we give a 
Weierstrass representation for these surfaces.  Proofs of the 
statements in this section may be found in \cite{HKR}.  

	For any Riemann example $\cal R$, there exists a translation $T$ of 
$\bfR^3$ so that $T$: leaves $\cal R$ invariant; is orientation-preserving on 
$\cal R$; generates the full cyclic orientation-preserving 
translation symmetry group of 
$\cal R$.  The quotient-surface ${\cal R}/T$ is a twice-punctured 
rectangular torus.  Thus, 
Riemann examples can be parametrized via a Weierstrass 
representation on such a torus.  
Fix $\lambda > 0$.  We first define $\bar{M_\lambda}$: 
\begin{equation} 
\bar{M_\lambda} = \{ \; (z,w) \in (\bfC \cup \{\infty\})^2 : 
w^2=z(z-\lambda)(z+\frac{1}{\lambda}) \; \} \; . 
\end{equation} 
Let $M_\lambda$ be the twice-punctured torus 
\begin{equation}
M_\lambda = \bar{M_\lambda} \setminus \{ (0,0), \, (\infty,\infty) \} \; . 
\end{equation}
With the Weierstrass data 
\begin{equation}
g(z,w) = z \; , \; \eta(z,w) = \frac{dz}{zw} \; , 
\end{equation}  
the Weierstrass representation for the Riemann examples 
${\cal R}_\lambda$ is 
\begin{equation}
R_\lambda(p) := Re\int_{p_0}^{p} \Phi_{(g,\eta)} \, dz \; , \; \; \; 
\Phi_{(g,\eta)} := \left( \begin{array}{c}
	\frac{(1-z^2)dz}{zw} \\
	\frac{i(1+z^2)dz}{zw} \\
	\frac{2dz}{w}
	\end{array}
\right) \; , \; \; \; p \in M_\lambda \; . 
\end{equation}
The family ${\cal R}_\lambda$ of Riemann examples is parametrized by 
$\lambda \in (0,\infty)$.  

	The surface $M_\lambda$ can be considered to be a double covering 
of the punctured $z$-plane, branched at $\lambda$ and $\frac{-1}{\lambda}$.  
Consider a circle $\alpha$ in the $z$-plane with center 
$\frac{-1}{2\lambda}$ and radius $\frac{1}{2}(\lambda + \frac{1}{\lambda})$, 
and lift $\alpha$ to a closed curve $\hat{\alpha}$ in $M_\lambda$.  
The Weierstrass integral about $\hat{\alpha}$ results in a 
nonzero period-vector, and this is the translation $T$.  

	Letting $\acute{z} = -z$ and $\acute{w} = iw$, we have 
\[ \acute{w}^2 = \acute{z}(\acute{z} - \frac{1}{\lambda}) 
(\acute{z} + \lambda) \; . \] 
Thus we may consider the map $(z,w) \rightarrow (\acute{z},
\acute{w})$ to be a conformal diffeomorphism between $M_\lambda$ and 
$M_{\frac{1}{\lambda}}$.  

	Because the metric on ${\cal R}_\lambda$ (\cite{Os} p65, 
\cite{HKR}) is 
\[ ds^2 = \frac{(1 + |z|^2)^2}{4|z|^2|w|^2} |dz|^2 \; , \] 
the maps $(z,w) \rightarrow (\bar{z},\bar{w})$ and 
$(z,w) \rightarrow (\bar{z},-\bar{w})$ both restrict to isometries of 
$M_\lambda$ in the induced metric.  The set 
${\cal A} = \{ (t,w) \in M_\lambda \, | \, t \in \bfR \}$, which consists of 
the union of their fixed-point sets, is mapped by the Weierstrass integral to 
geodesics in ${\cal R}_\lambda$.  It is clear from the 
Weierstrass integral that the intervals in $\cal A$ of the form 
$t \leq \frac{-1}{\lambda}$ or $0 < t \leq \lambda$ represent lines in 
${\cal R}_\lambda$, 
while intervals in $\cal A$ of the form $\frac{-1}{\lambda} 
\leq t < 0$ or $t \geq \lambda$ represent planar geodesics in 
${\cal R}_\lambda$.  Planar geodesics must be principal curves.  

	The surfaces ${\cal R}_\lambda$ 
have the following additional properties: 
\newcounter{num}
\begin{list}%
{\arabic{num})}{\usecounter{num}\setlength{\rightmargin}{\leftmargin}}

\item ${\cal R}_\lambda$ is 
foliated by circles and lines in horizontal planes.  
All the lines in ${\cal R}_\lambda$ are the image 
under the Weierstrass representation of points in $\cal A$ such 
that $t \leq \frac{-1}{\lambda}$ or $0 < t \leq \lambda$, which is the 
fixed-point set of the map $(z,w(z)) \rightarrow (\bar{z},-\bar{w}(z))$.  

\item ${\cal R}_\lambda$ has an infinite number of equally spaced horizontal 
flat ends.  These ends correspond to the punctures of $\bar{M_\lambda}$ at 
$z = 0, \infty$.  

\item ${\cal R}_\lambda$ is invariant under 
reflection in a plane parallel to the 
$(x_1,x_3)$-plane.  The intersection of this plane with ${\cal R}_\lambda$ 
consists of the image 
under the Weierstrass representation of points in $\cal A$ such 
that $\frac{-1}{\lambda} \leq z < 0$ or $z \geq \lambda$, which is the 
fixed-point set of the map $(z,w(z)) \rightarrow (\bar{z},\bar{w}(z))$.  

\item ${\cal R}_\lambda$ is 
invariant under rotation about horizontal lines that are 
perpendicular to the $(x_1,x_3)$-plane, and meet the surface orthogonally at 
the points where $g(z) = \pm i$.  
The map $(z,w) \rightarrow (\frac{-1}{z},\frac{w}{z^2})$ restricted to 
$M_\lambda$ represents this rotation.  

\end{list}

We will consider a Riemann example to be any surface that is the image of one 
of these surfaces ${\cal R}_\lambda$ 
under an isometry or homothety of $\bfR^3$.  

\subsection{Conjugate Pairs of Riemann Examples}

\begin{proposition}
	For any positive $\lambda \neq 1$, the Riemann examples 
${\cal R}_\lambda$ and 
${\cal R}_{\frac{1}{\lambda}}$ are conjugate but not congruent.  
${\cal R}_1$ is self-conjugate.
\end{proposition}

\begin{proof}
	Choose $\lambda > 1$.  The conjugate surface of ${\cal R}_\lambda$ has 
Weierstrass data 
\[g = z \; , \] \[i\eta = \frac{idz}{zw} = \frac{dz}{(-z)(iw)} \; . \] 
The Weierstrass data for the conjugate surface, 
expressed in terms of $(\acute{z},\acute{w}) = (-z,iw)$ on 
$M_{\frac{1}{\lambda}}$, is \[\tilde{g}(\acute{z},\acute{w}) = -\acute{z} 
\; , \] \[\tilde{\eta}(\acute{z},\acute{w}) = 
\frac{-d\acute{z}}{\acute{z}\acute{w}} \; . \] 
Because 
\[ \Phi_{(\tilde{g},\tilde{\eta})} = \left( \begin{array}{rrr}
	-1 & 0 & 0 \\
	0 & -1 & 0 \\
	0 & 0 & 1 
	\end{array}
\right) \Phi_{(g,\eta)} \; , \] 
it follows from equation (2.4) that 
the minimal surface in $\bfR^3$ with this Weierstrass data is the image 
of ${\cal R}_{\frac{1}{\lambda}}$ under a rotation by $\pi$ about a 
vertical axis.  Hence 
${\cal R}_\lambda$ and ${\cal R}_{\frac{1}{\lambda}}$ are conjugate surfaces.  

We use the Weierstrass data to show that ${\cal R}_\lambda$ and 
${\cal R}_{\frac{1}{\lambda}}$ cannot be congruent for $\lambda \neq 1$.  
Consider the points on 
${\cal R}_\lambda$ where $g = \pm 1$.  As shown in the previous
Section~2.1 (Property~1), 
when $\lambda > 1$, these points lie 
on the lines in ${\cal R}_\lambda$.  By contrast, the points of 
${\cal R}_{\frac{1}{\lambda}}$ where $g(z) = \pm 1$ lie on the planar 
geodesics in the vertical plane of reflective symmetry of the surface.  
Suppose there exists a symmetry of $\bfR^3$ taking ${\cal R}_\lambda$ to 
${\cal R}_{\frac{1}{\lambda}}$, 
then this symmetry must take lines to lines, ends to ends, and planar 
geodesics to planar geodesics.  
Since the lines in these two surfaces are all parallel to the 
$x_2$-axis and the ends are horizontal, the symmetry of $\bfR^3$ must 
be of the form 
\[ \left( \begin{array}{rrr}
\pm 1 & 0 & 0 \\
0 & \pm 1 & 0 \\
0 & 0 & \pm 1
\end{array} \right) \; . \] 
Therefore points of ${\cal R}_\lambda$ where 
$g = \pm 1$ must be mapped to points of ${\cal R}_{\frac{1}{\lambda}}$ where 
$g = \pm 1$.  This implies that the points of ${\cal R}_\lambda$ and 
${\cal R}_{\frac{1}{\lambda}}$ where $g = \pm 1$ are the 
points where the lines 
and planar geodesics in the surface intersect i.e.\ where $ g=\pm
\lambda ^{\pm 1} $.  But we are assuming $ \lambda \neq 1 $. 
Thus ${\cal R}_\lambda$ and 
${\cal R}_{\frac{1}{\lambda}}$ are not congruent.  
\end{proof}

\begin{remark}
One checks that the mapping $ (z,w)\rightarrow (z',w') $ from 
$ M_\lambda $ to $ M_{\frac{1}{\lambda}} $ has the property that 
on fundamental cycles it preserves the periods (one of which is
always zero) of (2.4).  Hence this mapping induces an isometry
of $ {\cal R}_\lambda $ and $ {\cal R}_{\frac{1}{\lambda}} $.
\end{remark}

Upon considering the behavior of ${\cal R}_\lambda$ and 
${\cal R}_{\frac{1}{\lambda}}$ as $\lambda \rightarrow 0$, 
one finds that ${\cal R}_\lambda$ and ${\cal R}_{\frac{1}{\lambda}}$ do 
not have limits that are regular surfaces.  However, if these surfaces are 
normalized 
appropriately they converge to a catenoid and its conjugate surface, the 
helicoid, which is consistent with Proposition~2.1.  
We shall prove this in the next section.  

\section{Proof of the Main Result}

To find a sequence of Riemann examples that converges to a catenoid or 
helicoid, we renormalize the surfaces ${\cal R}_\lambda$ as follows.  
Begin the integration in the Weierstrass representation 
at the point $z_0
 = 1$.  If $\lambda > 1$, rescale by 
$\sqrt{\lambda}$; if $\lambda < 1$, rescale 
by $\frac{1}{\sqrt{\lambda}}$.  We name 
these normalized Riemann examples $\acute{\cal R}_\lambda$.  Note that the 
rescaling is accomplished by multiplying $\eta$ in the Weierstrass 
representation by $\sqrt{\lambda}$ and $\frac{1}{\sqrt{\lambda}}$, 
respectively.  Thus, the 
Weierstrass representation for $\acute{\cal R}_\lambda$ is 
\[ g = z \; , \; \eta = \frac{\sqrt{\lambda}dz}{zw} \; , \; \lambda \geq 1 \; 
, \] 
\[ g = z \; , \; \eta = \frac{dz}{\sqrt{\lambda}zw} \; , \; \lambda \leq 1 \; 
. \]  One can check as in Proposition~2.1 that 
$\acute{\cal R}_\lambda$ and 
$\acute{\cal R}_{\frac{1}{\lambda}}$ are conjugate. 

Let $ K $ denote Gauss curvature.  


\begin{lemma}
	There exists a universal bound $c$ such that $|K| \leq c$ on 
$\acute{\cal R}_\lambda$ for all $\lambda \in (0,\infty)$.  
\end{lemma}

\begin{proof}
With Weierstrass data given locally as $g$ and $\eta = fdz$, 
the Gauss curvature of a minimal surface is (\cite{Os}, p76) 
\begin{equation} K = -(\frac{4|g^\prime|}{|f|(1+|g|^2)^2})^2 \; . 
\end{equation}
In the case of $\acute{\cal R}_\lambda$ a computation using equations (2.1) 
and (3.5) yields 
\[|K| = 16\lambda^{\alpha}\frac{|z-\lambda||z+{\lambda^{-1}}|}{
|z|(|z|+|z|^{-1})^4} \; . \] 
where $\alpha = 1$ if $\lambda \leq 1$, and $\alpha = -1$ if $\lambda 
\geq 1$.  If $\lambda \leq 1$, then 
\[|K| \leq 16\frac{(|z|+1)\lambda(|z|+{\lambda^{-1}})}{
|z|(|z|+|z|^{-1})^4} \; \leq \; 
16\frac{(|z|+1)^2}{|z|(|z|+|z|^{-1})^4} \; \leq 4 \; . \] 
Similarly, if $\lambda \geq 1$, then 
\[|K| \leq 16\frac{{\lambda^{-1}}(|z|+\lambda)(|z|+1)}{
|z|(|z|+|z|^{-1})^4} \; \leq \; 
16\frac{(|z|+1)^2}{|z|(|z|+|z|^{-1})^4} \; \leq 4 \; . \] 
Thus, $|K|$ has a universal upper bound.  
\end{proof}

\begin{remark}
More calculation shows that $|K|$ on ${\cal R}_\lambda$ 
is maximized at a fixed point of one of the 
symmetries of the surface.  Thus $K$ is maximized in one of the 
following three places: somewhere along a straight line in 
${\cal R}_\lambda$; somewhere along a planar geodesic in 
the plane of reflective symmetry parallel to $\{x_2 = 0\}$; or at the 
points where $g(z) = \pm i$, which are the fixed points of the normal 
rotation.  At the points where $z = g(z) = \pm i$, the value of 
$|K|$ on ${\cal R}_\lambda$ is $\lambda + \frac{1}{\lambda}$.  Hence
on $ \acute{\cal R}_\lambda $ the value is $ 1+\frac{1}{\lambda^{2}}
(\mbox{resp. }1+\lambda^2) $ when $ \lambda > 1 (\mbox{resp.\ when }
\lambda < 1) $.  

The following conjecture has been verified numerically:  
The absolute value of the Gaussian curvature on ${\cal R}_\lambda$ 
is maximized at the points where $z = g(z) = \pm i$.  This implies
that the optimal value in Lemma~3.1 is $ c=2 $.
\end{remark}

\begin{figure}[p]
   \hspace{2in}
   \epsfxsize = 1.7in	\epsffile{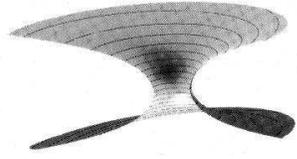}
\hfill \break
   \caption{Portion of the Riemann Example 
	$\acute{\cal R}_\lambda, \; \lambda = 0.1$}
\end{figure}

\begin{figure}[p]
   \hspace{2in}
   \epsfxsize = 1.7in	\epsffile{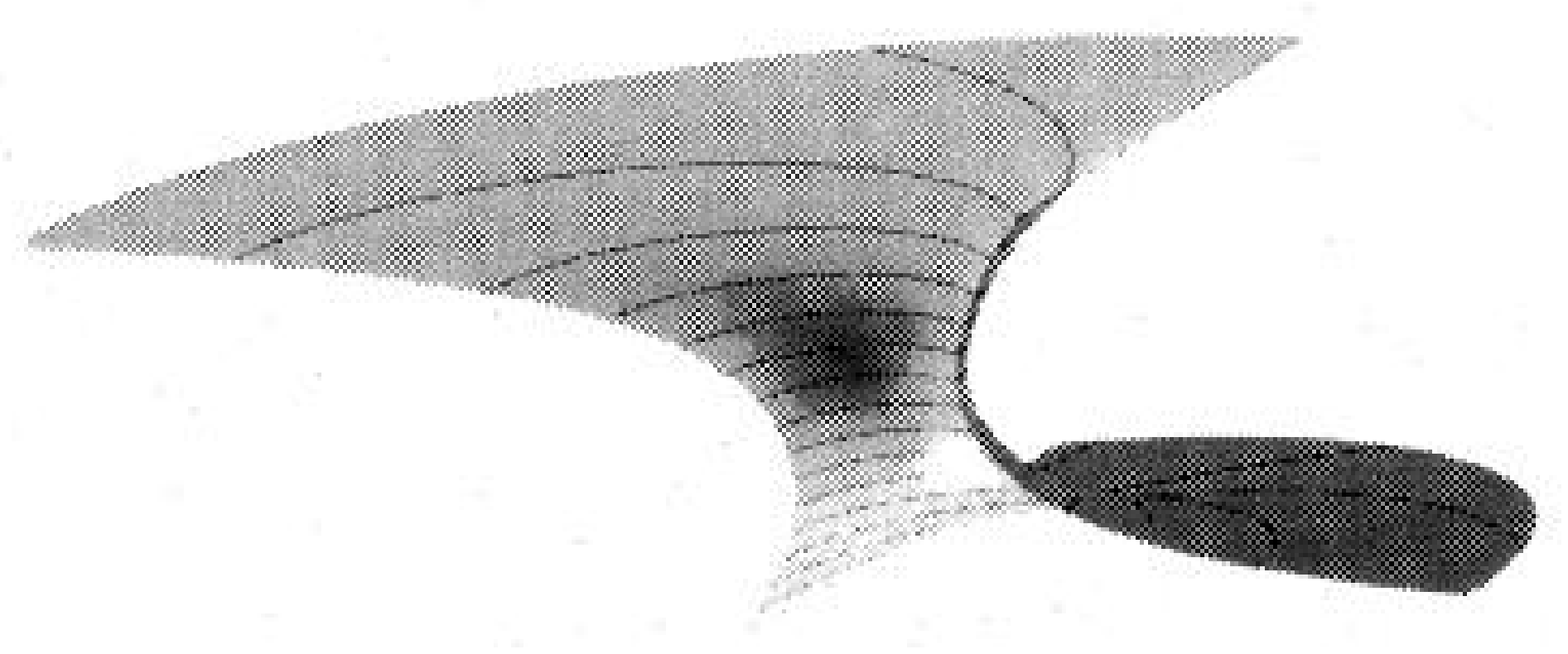}
\hfill \break
   \caption{Portion of the Riemann Example 
	$\acute{\cal R}_\lambda, \; \lambda = 0.5$}
\end{figure}

\begin{figure}[p]
   \hspace{2in}
   \epsfxsize = 1.7in	\epsffile{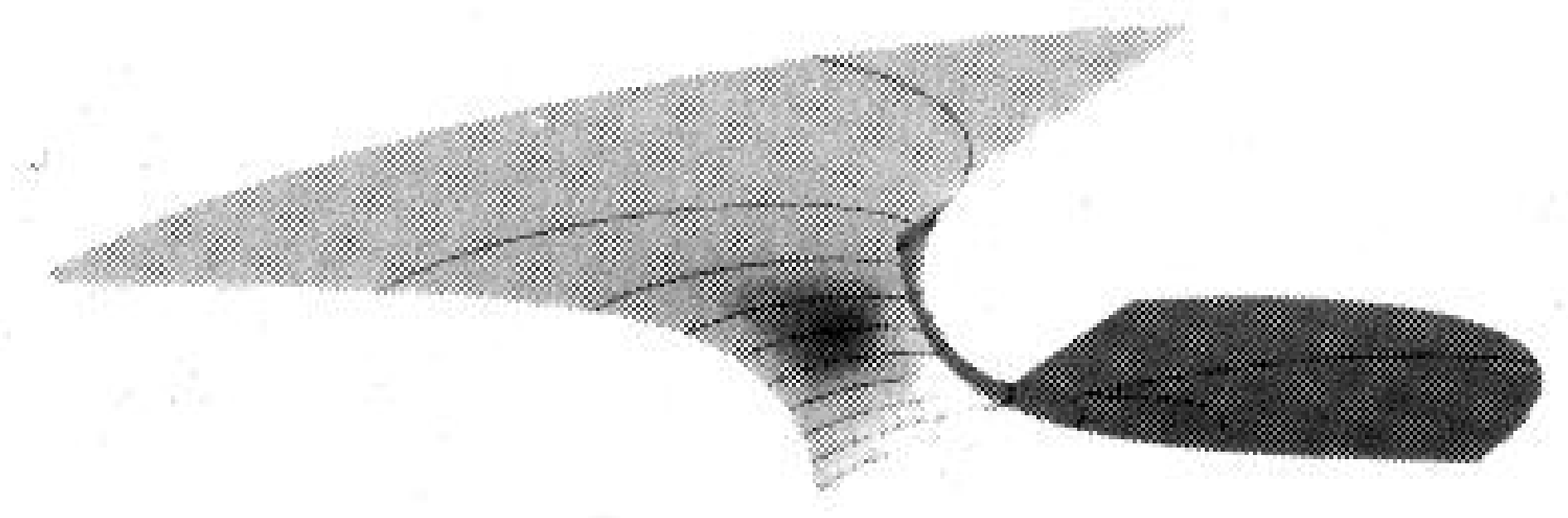}
\hfill \break
   \caption{Portion of the Riemann Example 
	$\acute{\cal R}_\lambda, \; \lambda = 1.0$}
\end{figure}

\begin{figure}[p]
   \hspace{2in}
   \epsfxsize = 1.7in	\epsffile{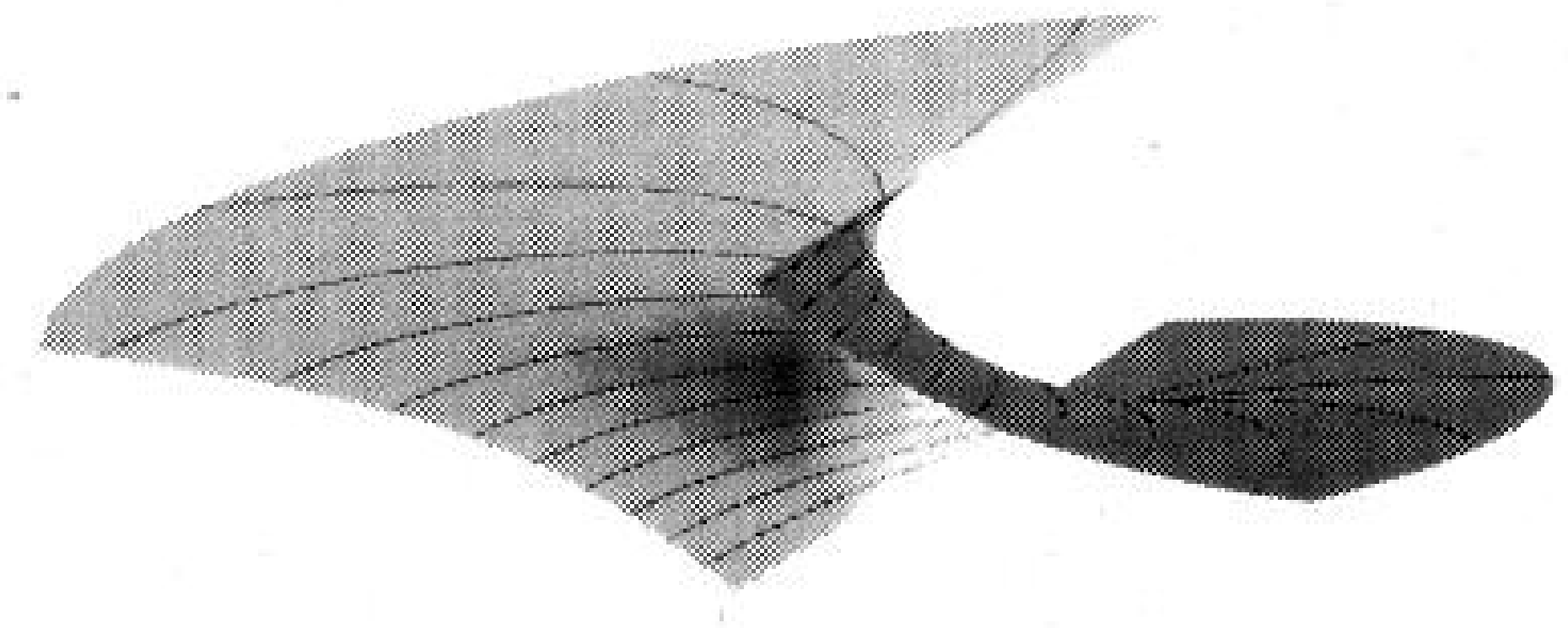}
\hfill \break
   \caption{Portion of the Riemann Example 
	$\acute{\cal R}_\lambda, \; \lambda = 3.0$}
\end{figure}

\begin{figure}[p]
   \hspace{2in}
   \epsfxsize = 1.7in	\epsffile{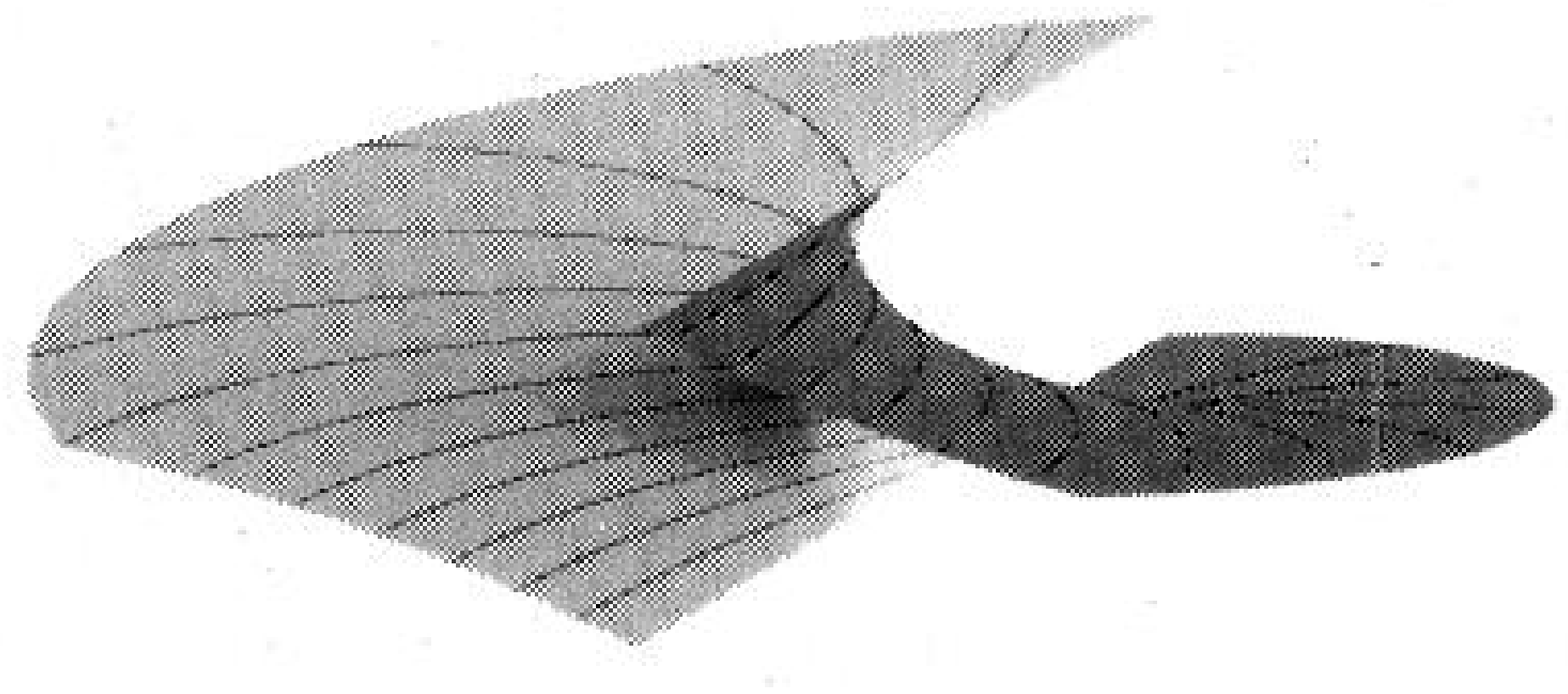}
\hfill \break
   \caption{Portion of the Riemann Example 
	$\acute{\cal R}_\lambda, \; \lambda = 5.0$}
\end{figure}

\begin{figure}[p]
   \hspace{2in}
   \epsfxsize = 1.7in	\epsffile{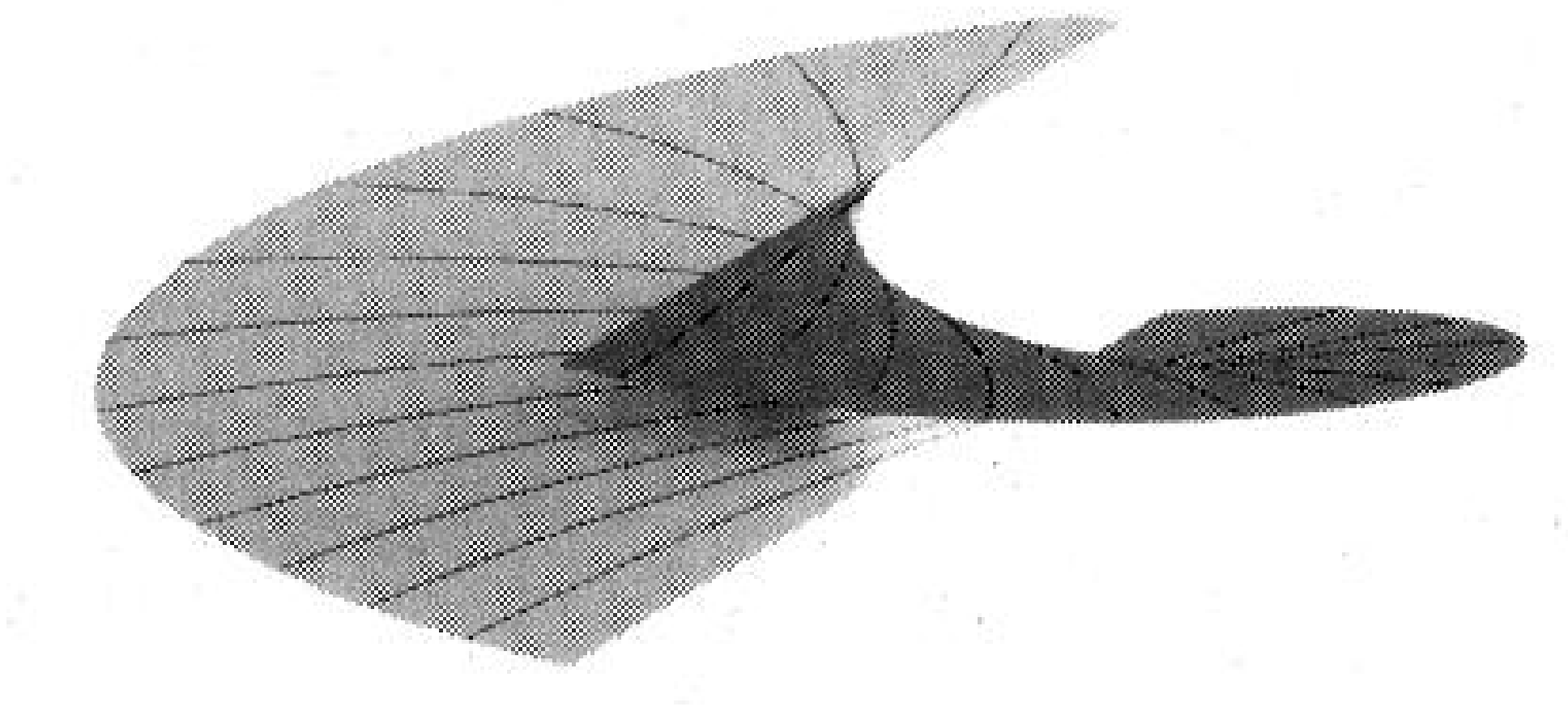}
\hfill \break
   \caption{Portion of the Riemann Example 
	$\acute{\cal R}_\lambda, \; \lambda = 10.0$}
\end{figure}

\begin{proposition}
	There exist sequences of Riemann examples that converge to a 
helicoid and sequences that converge to a catenoid.
\end{proposition}

\begin{proof}
	We will prove this proposition by showing that the 
limit, as $\lambda \rightarrow 0$, of $\acute{\cal R}_\lambda$ is a 
catenoid, and the limit, as $\lambda \rightarrow \infty$, of 
$\acute{\cal R}_\lambda$ is a helicoid (see Figures 2 through 4).  We 
consider first the case $\lambda \rightarrow 0$.  

	For a large positive number 
$L$, let $A_L$ be the annular ring 
$\{z \in \bfC : \frac{1}{L} < |z| < L\}$.  One can think of $z$ in
(2.1) as a map 
from $M_\lambda$ to $\bfC \setminus \{0\}$.  Then 
for $\lambda$ sufficiently 
small, $z^{-1}(A_L)$ consists of two disjoint sets 
in $M_\lambda$.  Choose either one and name it $\acute{A}_{L,\lambda}$.  
Since $z$ parametrizes $\acute{A}_{L,\lambda}$, we shall refer to points 
in $\acute{A}_{L,\lambda}$ by their $z$ coordinates.  

Let $B_r$ be a ball centered at the origin in $\bfR^3$ with radius $r$.  
Let $S_{r}$ be the horizontal slab 
$\{(x_1,x_2,x_3) \in \bfR^3 \, | \, -r \leq x_3 \leq r\}$, 
with boundary planes $x_3 = \pm r$.  

Let $\cal C$ be the catenoid obtained by using the Weierstrass data 
\[ g(z) = z \; , \; \eta(z) = \frac{dz}{z^2} \; , \] 
and integrating over $\bfC \setminus \{0\}$ with base point $z_0 = 1$.  
The Weierstrass representation for the catenoid $\cal C$ is 
\[C(p) := Re\int_{1}^{p} \left( \begin{array}{c}
	\frac{(1-z^2)dz}{z^2} \\
	\frac{i(1+z^2)dz}{z^2} \\
	\frac{2dz}{z}
	\end{array}
\right) \; , \; p \in \bfC \setminus \{0\} \; .\]
Let $\hat{\cal C}_L$ be the image of the restriction of $C(p)$ 
to $A_L$.  By examination of the 
third coordinate of the Weierstrass representation, we see that 
$\hat{\cal C}_L = {\cal C} \cap S_{2 \ln L}$.  
Choose a large value for $r$, and 
choose $L > e^{r}$  so that ${\cal C} \cap S_{2r} \subseteq \hat{\cal C}_L$.  

	For $\lambda < 1$, the Weierstrass representation 
for $\acute{\cal R}_\lambda$ is 
\begin{equation}
\acute{R}_\lambda(p) := Re\int_{1}^{p} \left( \begin{array}{c}
	\frac{(1-z^2)dz}{\sqrt{\lambda}zw} \\
	\frac{i(1+z^2)dz}{\sqrt{\lambda}zw} \\
	\frac{2dz}{\sqrt{\lambda}w}
	\end{array}
\right) = C(p) + 
Re\int_{1}^{p} \left( \begin{array}{c}
	f_{0,\lambda}(z)\frac{(1-z^2)dz}{z^2} \\
	f_{0,\lambda}(z)\frac{i(1+z^2)dz}{z^2} \\
	f_{0,\lambda}(z)\frac{2dz}{z}
	\end{array}
\right) \; , 
\end{equation} 
where 
\[ f_{0,\lambda}(z) = \frac{\sqrt{z}}{\sqrt{(z-\lambda)(\lambda z+1)}}
 - 1 \; . \]

	If $\lambda$ is close to zero, then the Weierstrass integral 
$\acute{R}_\lambda$ along any closed curve in 
$\acute{A}_{L,\lambda}$ is zero.  Therefore the Weierstrass 
integral depends only on the endpoints of a path 
in $\acute{A}_{L,\lambda}$.  All paths from $z_0 = 1$ to $p \in A_L$ may be 
chosen to be of length less than $cL$, where $c$ is some fixed constant.  

	Note that the functions 
\[ \frac{1-z^2}{z^2} \; , \; \frac{1+z^2}{z^2} \; , \; 
\frac{2}{z} \] are bounded on $A_L$.  
Also, for all $z \in A_L$, 
\[\lim_{\lambda \rightarrow 0}|f_{0,\lambda}(z)| = 0 \; . \]  
It follows that the last integral in 
equation~(3.6) becomes arbitrarily small, 
uniformly on $\acute{A}_{L,\lambda}$, as $\lambda \rightarrow 0$.  
In particular, for all $p \in A_L$, $\lim_{\lambda \rightarrow 0} 
\acute{R}_\lambda(p) = C(p)$.  

Let $\hat{\cal R}_\lambda$ be the image of $\acute{A}_{L,\lambda}$ 
under the map $\acute{R}_\lambda(p)$.  It follows that 
$\hat{\cal R}_\lambda$ is a small perturbation of $\hat{\cal C}_L$ 
for $\lambda$ close to zero.  We conclude that 
$\acute{\cal R}_\lambda \cap B_r \subseteq \hat{\cal R}_\lambda$.  
Thus to examine the behavior of $\acute{\cal R}_\lambda$, 
as $\lambda \rightarrow 0$, in arbitrary compact 
regions of $\bfR^3$, it is enough to examine 
$\hat{\cal R}_\lambda$ in $B_r$ for arbitrarily large fixed $r$.  

We claim that the surface $\acute{\cal R}_\lambda \cap B_r$ is a 
graph over a subset of the catenoid $\cal C$, for $\lambda$ close to $0$.  
By Lemma~3.1, the Gaussian curvature $K$ is uniformly bounded for 
all $\acute{\cal R}_\lambda$.  This implies that the normal curvatures are
also uniformly founded.  
This in turn implies $\acute{\cal R}_\lambda \cap B_r$ is a 
union of graphs over a subset of the catenoid $\cal C$.  
Since the surface $\acute{\cal R}_\lambda$ is foliated by 
lines and circles in horizontal planes, there is a single graph.  

We now show that $\acute{\cal R}_\lambda \cap B_r$ converges to 
${\cal C} \cap B_r$, in the $C^\infty$-topology, as 
$\lambda \rightarrow 0$.  Because the curvature of $\acute{\cal R}_\lambda $ 
is bounded uniformly in $ \lambda $, it follows from standard
elliptic theory that the convergence is $ C^\infty $.  (Curvature
bounds give $ C^1 $ estimates.)  However, we can give in our
specific case a direct proof.  Let $p \in {\cal C} \cap B_r$.  
Let $z_p \in A_L$ be the point such that $C(z_p) = p$, and let 
$B_h(z_p) \subseteq A_L$ be a ball of radius $h$ about $z_p$.  Note that 
$\acute{R}_\lambda(B_h(z_p))$ 
consists of disks in $\bfR^3$ converging in the $C^0$-norm, 
as $\lambda \rightarrow 0$, 
to a disk containing $p$ and lying on the catenoid $\cal C$.  Let 
$x_{i,\lambda}(z)$ be the $i$'th coordinate function of 
$\acute{\cal R}_\lambda$, let $x_i(z)$ be the $i$'th coordinate 
function on $\cal C$, 
and let $f_\lambda(z) = x_{i,\lambda}(z) - x_i(z)$.  Note that 
\[ \max_{w \in B_h(z_p)} \left| f_\lambda(w) \right| \rightarrow 0 \, ,
\mbox{ as } \lambda \rightarrow 0 \; . \] 
Since $f_\lambda(z)$ is harmonic, it is the real part of a holomorphic 
function $F_\lambda(z)$, and $F_\lambda(z)$ can be chosen so that 
\[ \max_{w \in B_h(z_p)} \left| F_\lambda(w) \right| \rightarrow 0 \, ,
\mbox{ as } \lambda \rightarrow 0 \; . \]  
Assuming $z \in B_{\frac{h}{2}}(z_p)$, and using the Cauchy integral formula, 
we have 
\[ \left| f^{(k)}_\lambda(z) \right| \leq \left| F^{(k)}_\lambda(z) \right| 
= \left| \frac{k!}{2\pi i} \int_{\partial B_h(z_p)} \frac{F(\xi)}
{(\xi - z)^{k+1}} d\xi \right| \] 
\[ \leq \frac{k!}{2\pi} \left(\frac{2}{h}\right)^{k+1} 2\pi h 
\max_{w \in B_h(z_p)} \left| F_\lambda(w) \right| \; , \] 
and so we have that $|f^{(k)}_\lambda(z)| \rightarrow 0 \; \forall k, \; 
\forall z \in B_{\frac{h}{2}}(z_p)$, 
as $\lambda \rightarrow 0$.  Therefore $|f^{(k)}_\lambda(z)| 
\rightarrow 0 \; \forall k, \; \forall z \in A_L$, 
and the convergence is in the $C^\infty$-topology in $B_r$.

Since $r$ is arbitrary, we have completed the proof that 
$\acute{\cal R}_\lambda$ 
converges to the catenoid $\cal C$ as $\lambda \rightarrow 0$.  

	It is now intuitively clear that the surfaces 
$\acute{\cal R}_\lambda$ converge to 
a helicoid for the following reason: as $\lambda \rightarrow \infty$, 
their conjugate surfaces $\acute{\cal R}_{\frac{1}{\lambda}}$ 
converge to a catenoid.  One can give an explicit proof of this, 
similar to the proof just given for the case $\lambda \rightarrow 
0$.  However, there are some differences in the proof of the 
case $\lambda \rightarrow \infty$.    
For $\lambda$ close to $\infty$, the representation for $\acute{\cal R}_
\lambda$ is 
\[ \acute{R}_\lambda(p) := Re\int_{1}^{p} \left( \begin{array}{c}
	\frac{\sqrt{\lambda}(1-z^2)dz}{zw} \\
	\frac{i\sqrt{\lambda}(1+z^2)dz}{zw} \\
	\frac{2\sqrt{\lambda}dz}{w}
	\end{array}
\right) = \]
\begin{equation}
-Re\int_{1}^{p} \left( \begin{array}{c}
	\frac{i(1-z^2)dz}{z^2} \\
	\frac{-(1+z^2)dz}{z^2} \\
	\frac{2idz}{z}
	\end{array}
\right) + 
 Re\int_{1}^{p} \left( \begin{array}{c}
	f_{\infty,\lambda}(z)\frac{i(1-z^2)dz}{z^2} \\
	f_{\infty,\lambda}(z)\frac{-(1+z^2)dz}{z^2} \\
	f_{\infty,\lambda}(z)\frac{2idz}{z}
	\end{array}
\right) \; ,  
\end{equation}  
where 
\[ f_{\infty,\lambda}(z) = 1 - 
\frac{\sqrt{z}}{\sqrt{(1-\frac{z}{\lambda})(z+\frac{1}{\lambda})}}
 \; .  \] 

Note that the first integral in the above sum in 
equation~(3.7) is the Weierstrass integral for a helicoid, 
since Weierstrass data for a helicoid is 
\[ g(z) = z \; , \; \eta(z) = \frac{idz}{z^2} \; . \]

If $\lambda$ is 
large, then the vertical distance between adjacent planar ends of 
$\acute{\cal R}_\lambda$ is approximately $2\pi$.  This can be verified by 
integrating the third function in the integrand 
of the Weierstrass representation from $+1$ to 
$-1$ along the half of the unit circle lying 
in the upper half of the complex plane.  
This curve is $z = e^{it}, 0 \leq t \leq \pi$.  The distance between adjacent 
planar ends for large $\lambda$ is 
\[Re\int_0^\pi\frac{2\sqrt{\lambda}ie^{it}dt}
{\sqrt{e^{it}(e^{it}-\lambda)(e^{it}+\frac{1}{\lambda})}} \approx 
Re\int_0^\pi\frac{2\sqrt{\lambda}idt}{\sqrt{-\lambda}} = 2\pi \; . \]

The major difference from the proof in the case $\lambda 
\rightarrow 0$ is this: When $\lambda$ is close to 
$\infty$, then a homologically nontrivial loop in $\acute{A}_
{L,\lambda}$ has a nonzero real period with respect to the Weierstrass 
integral.  Now one may only assume that a path from $z=1$ to 
$z=p$ in $A_L$ has length less than 
$c n L$, where $n$ is the number of 
times that the path wraps about the origin, and $c$ is some fixed constant.  
Since the distance between adjacent ends is approximately $2\pi$, 
it follows that if the wrapping number $n$ about the origin 
of a path from 1 to $p$ is sufficiently large, then $z=p$ will be mapped by 
the Weierstrass integral to a point outside $B_r$.  
Thus, we may assume there is an upper bound $N$ for $|n|$, depending only 
on $r$.  Hence, for large $\lambda$, we may assume that any path from 1 
to $p$ has length less than $\acute{c} L$, where $\acute{c} = c N$.  

The other parts of the proof of the case $\lambda \rightarrow 
\infty$ transfer directly from the proof of the case $\lambda \rightarrow 0$.  
\end{proof}

\begin{corollary}
	As $\lambda \rightarrow \infty$, the radius of any  
level-circle of $\acute{\cal R}_\lambda$ diverges to infinity.  
\end{corollary}

\begin{corollary}
	The planar 
curve described by the set of centers of the horizontal 
circles that foliate $\acute{\cal R}_\lambda$ has the property that
its maximum curvature  
approaches zero as $\lambda \rightarrow \infty$.  
\end{corollary}

	It remains to determine what other surfaces are limits 
of Riemann examples.  Three possibilities are clear: 
a single flat plane; an infinite number of equally-spaced flat planes; and 
another Riemann example.  The following proposition states that these are 
the only other possibilities, completing the proof of Theorem~1.1.  

\begin{proposition}
Any convergent sequence of Riemann examples converges to one of  
the following surfaces:

	a Riemann example;

	a helicoid;

	a catenoid;

	a single flat plane;

	an infinite number of equally-spaced flat planes.
\end{proposition}

\begin{proof}
	Let $\{{\cal R}_j\}_{j=1}^\infty$ be a sequence of Riemann examples 
converging to a surface $\cal S$.  From the definition of convergence, 
$\cal S$ is clearly minimal, properly embedded, and complete.  
Without loss of generality we may assume that the
limiting normals $\vec{v_j}$ at the 
ends of the surfaces ${\cal R}_j$ converge in $\bf S^2$ to
a vertical vector.  With this normalization, the surfaces ${\cal R}_j$ are 
foliated by circles and lines lying in planes that are becoming 
horizontal as $j \rightarrow \infty$.  These almost-horizontal planes 
intersect ${\cal R}_j$ 
in connected sets, and these connected sets each consist of a 
single circle or a single line.  Thus every curve ${\cal S} \cap \{x_3 = c\}$ 
is the limit set of a sequence of 
circles and/or lines.  It follows that $ {\cal S} \cap \{x_{3}=c\} $ is
either a circle, a line, the empty set or the plane $ \{x_{3}=c\} $.
Recall from the introduction that this implies that each component of
$ \cal S $ is either a helicoid, a catenoid, a Riemann example or a
plane.  Moreover: the helicoid and catenoid must have vertical
axes; the plane and the ends of the Riemann example must be horizontal.
If $ \cal S $ has more than one component, they must all be horizontal
planes (because $ \cal S $ is fibred by circles or lines, one in each
horizontal plane).  Suppose $ \cal S $ consists of two or more horizontal
planes.
Since the points on the surfaces ${\cal R}_j$ where 
the Gauss map is vertical must 
approach the ends of ${\cal R}_j$ as $j \rightarrow \infty$, 
it is the equally spaced-ends of the surfaces ${\cal R}_j$ that 
converge to $\cal S$.  
Therefore $\cal S$ is an {\em infinite} collection of equally-spaced planes.  
\end{proof}

\end{document}